\documentclass[11pt]{article}
\usepackage{amsmath, amssymb, amsfonts, amsthm, latexsym} 
\usepackage{mathrsfs}
\usepackage{graphicx}
\usepackage{authblk}
\usepackage[usenames]{color}
\newtheorem{main}{Theorem}
\newtheorem{theorem}{Theorem}[section]
\newtheorem{lemma}[theorem]{Lemma}

\newtheorem{definition}[theorem]{Definition}

\newtheorem{remark}[theorem]{Remark}

\newtheorem*{notat*}{Notation}

\title{The Hamilton cycle space of random graphs}

\author{Dan Hefetz \thanks{School of Computer Science, Ariel University, Ariel 40700, Israel. Email: {\tt danhe@ariel.ac.il}.}
\quad Michael Krivelevich \thanks{School of Mathematical Sciences, Tel Aviv University, Tel Aviv 6997801, Israel. Research supported in part
by NSF-BSF grant 2023688. Email: {\tt krivelev@tauex.tau.ac.il}.}}

\begin{document}
\maketitle
 
\begin{abstract}
The cycle space of a graph $G$, denoted $\mathcal{C}(G)$, is a vector space over ${\mathbb F}_2$, spanned by all incidence vectors of edge-sets of cycles of $G$. If $G$ has $n$ vertices, then $\mathcal{C}_n(G)$ denotes the subspace of $\mathcal{C}(G)$, spanned by the incidence vectors of Hamilton cycles of $G$. A classical result in the theory of random graphs asserts that for $G \sim \mathbb{G}(n,p)$, asymptotically almost surely the necessary condition $\delta(G) \geq 2$ is also sufficient to ensure Hamiltonicity. Resolving a problem of Christoph, Nenadov, and Petrova, we augment this result by proving that $\mathcal{C}_n(G) = \mathcal{C}(G)$ holds asymptotically almost surely for $G \sim \mathbb{G}(n,p)$, with $n$ being odd and $p := p(n) \geq (\log n + 2 \log \log n + \omega(1))/n$. That is, not only does $G$ typically have a Hamilton cycle, but its Hamilton cycles are typically rich enough to span its cycle space. It is known and easy to show that if $p := p(n) \leq (\log n + 2 \log \log n - \omega(1))/n$, then a.a.s. $\mathcal{C}_n(G) \neq \mathcal{C}(G)$ unless $G$ is a forest.  
\end{abstract}


\section{Introduction}  
Let $G = (V,E)$ be a graph on $n$ vertices. The \emph{edge space} of $G$, denoted $\mathcal{E}(G)$, is a vector space over ${\mathbb F}_2$ consisting of all incidence vectors of subsets of $E$. The \emph{cycle space} of $G$, denoted $\mathcal{C}(G)$, is the subspace of $\mathcal{E}(G)$, spanned by all incidence vectors of cycles of $G$. For any integer $3 \leq k \leq n$, let $\mathcal{C}_k(G)$ be the subspace of $\mathcal{C}(G)$, spanned by all incidence vectors of cycles of length $k$ in $G$. Determining conditions under which $\mathcal{C}_k(G) = \mathcal{C}(G)$ holds for some $3 \leq k \leq n$ is a well-studied problem (see, e.g.,~\cite{BK, BL, DHJ, Hartman, Locke1, Locke2}). In this paper we are interested in the case $k = n$, that is, in graphs whose cycle space is spanned by their Hamilton cycles. This problem has been addressed by various researchers (see, e.g.,~\cite{ALW, Heinig1, Heinig2}). Since the symmetric difference of any two even graphs (i.e., subsets of $E$ of even size) is an even graph, it is evident that if $\mathcal{C}_n(G) = \mathcal{C}(G)$, then $G$ is bipartite or $n$ is odd. Moreover, $G$ must either be acyclic or Hamiltonian. Since the former case is not very interesting, this study can be viewed as part of the common theme of proving that (possibly slightly strengthened) various sufficient conditions for Hamiltonicity in fact ensure stronger properties. 

Random graphs provide a natural setting for this problem. Indeed, already Erd\H{o}s and R\'enyi~\cite{ER} raised the question of what the threshold probability of Hamiltonicity in random graphs is. After a series of efforts by various researchers, including Korshunov~\cite{Korshunov} and P\'osa~\cite{Posa}, the
problem was finally solved by Koml\'os and Szemer\'edi~\cite{KSz} and independently by Bollob\'as~\cite{Bollobas}, who proved that if $p \geq (\log n + \log \log n + \omega(1))/n$, then $\mathbb{G}(n,p)$ is asymptotically almost surely (a.a.s. hereafter) Hamiltonian. Note that this is best possible as if $p \leq (\log n + \log \log n - \omega(1))/n$, then a.a.s. there are vertices of degree at most one in $\mathbb{G}(n,p)$. That is, for random graphs, the clearly necessary condition of having minimum degree at least two is a.a.s. sufficient for Hamiltonicity. Since then, it has been shown that random graphs are robustly Hamiltonian in various ways (see, e.g.,~\cite{BFHK, DGMS, FKL, HKLO, KKO, KSa, LS, Mont}).

As for the cycle space of random graphs, it was proved by Heinig~\cite{Heinig2} that $\mathcal{C}_n(\mathbb{G}(n,p)) = \mathcal{C}(\mathbb{G}(n,p))$ holds a.a.s. provided $n$ is odd and $p \geq n^{- 1/2 + o(1)}$. This was significantly improved by Christoph, Nenadov, and Petrova~\cite{CNP} who proved that a.a.s. $p \geq C \log n/n$ suffices, where $C$ is a sufficiently large constant. One may wonder if, as for Hamiltonicity, $\delta(\mathbb{G}(n,p)) \geq 2$ is a.a.s. sufficient to ensure  $\mathcal{C}_n(\mathbb{G}(n,p)) = \mathcal{C}(\mathbb{G}(n,p))$. However, as observed by Heinig~\cite{Heinig2}, if $p := p(n)$ is large enough to ensure that $\mathbb{G}(n,p)$ is a.a.s. not a forest, then for $\mathcal{C}_n(\mathbb{G}(n,p)) = \mathcal{C}(\mathbb{G}(n,p))$ to hold a.a.s., $p$ needs to be at least large enough so as to ensure that a.a.s. $\delta(\mathbb{G}(n,p)) \geq 3$. It was then asked by Christoph, Nenadov, and Petrova in~\cite{CNP} whether $\delta(\mathbb{G}(n,p)) \geq 3$ is a.a.s. sufficient. Our main result shows that if $p$ is large enough to ensure that
$\delta(\mathbb{G}(n,p)) \geq 3$ holds a.a.s., then it also ensures that $\mathcal{C}_n(G) = \mathcal{C}(G)$ holds a.a.s. (provided that $n$ is odd).

\begin{main} \label{th::CycleSpaceGnp}
Let $G \sim \mathbb{G}(n,p)$, where $n$ is odd and $p := p(n) \geq (\log n + 2 \log \log n + \omega(1))/n$. Then, a.a.s. $\mathcal{C}_n(G) = \mathcal{C}(G)$. 
\end{main}

Note that Theorem~\ref{th::CycleSpaceGnp} is indeed best possible as, for $p := p(n) \leq (\log n + 2 \log \log n - \omega(1))/n$, a.a.s. $\delta(\mathbb{G}(n,p)) \leq 2$ and thus, as noted above, a.a.s. $\mathcal{C}_n(\mathbb{G}(n,p)) \neq \mathcal{C}(\mathbb{G}(n,p))$ (unless it is a forest). 

\bigskip

The rest of this paper is organised as follows. In Section~\ref{sec::prelim} we introduce some terminology, notation, and standard tools, and present the method of Christoph, Nenadov, and Petrova from~\cite{CNP} which is a central ingredient in our proof. In Section~\ref{sec::random} we prove our main result, namely Theorem~\ref{th::CycleSpaceGnp}.

\section{Preliminaries and tools} \label{sec::prelim}

For the sake of simplicity and clarity of presentation, we do not make a particular effort to optimize the constants obtained in some of our proofs. We also omit floor and ceiling signs whenever these are not crucial. Most of our results are asymptotic in nature and whenever necessary we assume that the number of vertices $n$ is sufficiently large. Throughout this paper, $\log$ stands for the natural logarithm, unless explicitly stated otherwise. Our graph-theoretic notation is standard; in particular, we use the following.

For a graph $G$, let $V(G)$ and $E(G)$ denote its sets of vertices and edges respectively, and let $v(G) = |V(G)|$ and $e(G) = |E(G)|$. For a set $A \subseteq V(G)$, let $E_G(A)$ denote the set of edges of $G$ with both endpoints in $A$, and let $e_G(A) = |E_G(A)|$. For disjoint sets $A, B \subseteq  V(G)$, let $E_G(A, B)$ denote the set of edges of $G$ with one endpoint in $A$ and one endpoint in $B$, let $e_G(A, B) = |E_G(A, B)|$, and let $G[A,B]$ denote the bipartite graph with vertex-set $A \cup B$ and edge-set $E_G(A,B)$. For a set $S \subseteq V(G)$, let $G[S]$ denote the subgraph of $G$ induced by the set $S$. For a set $S \subseteq V(G)$, let $N_G(S) = \{v \in V(G) \setminus S : \exists u \in S \textrm{ such that } uv \in E(G)\}$ denote the \emph{external neighbourhood} of $S$ in $G$. For a vertex $u \in V(G)$ we abbreviate $N_G(\{u\})$ as $N_G(u)$ and let $\textrm{deg}_G(u) = |N_G(u)|$ denote the degree of $u$ in $G$. The maximum degree of a graph $G$ is $\Delta(G) = \max\{\textrm{deg}_G(u) : u \in V(G)\}$, and the minimum degree of a graph $G$ is $\delta(G) = \min \{\textrm{deg}_G(u) : u \in V(G)\}$. For a vertex $u \in V(G)$ and a set $S \subseteq V(G)$, let $N_G(u, S) = N_G(u) \cap S$ and let $\textrm{deg}_G(u, S) = |N_G(u, S)|$. For a vertex $x \in V(G)$, let $\partial_G(x) = \{xy : y \in N_G(x)\}$. Given any two (not necessarily distinct) vertices $x, y \in V(G)$, the \emph{distance} between $x$ and $y$ in $G$, denoted $\textrm{dist}_G(x,y)$, is the length of a shortest path between $x$ and $y$ in $G$, where the length of a path is the number of its edges (for the sake of formality, we define $\textrm{dist}_G(x,y)$ to be $\infty$ whenever $x$ and $y$ lie in different connected components of $G$). The \emph{diameter} of $G$, denoted $\textrm{diam}(G)$, is $\max \{\textrm{dist}_G(x,y) : x, y \in V(G)\}$. 



Throughout this paper we make use of the standard upper bound $\binom{n}{k} \leq \left(\frac{e n}{k} \right)^k$, holding for every integer $1 \leq k \leq n$, and of the following Chernoff-type bounds (see, e.g., ~\cite{JLR}).
\begin{theorem} \label{th::Chernoff}
Let $X \sim \emph{Bin}(n,p)$ and let $0 < \alpha < 1 < \beta$. Then,
\begin{enumerate}
\item [$(a)$] $\mathbb{P}(X \leq \alpha \mathbb{E}(X)) \leq \exp \{- (\alpha \log \alpha - \alpha + 1) \mathbb{E}(X)\}$, and

\item [$(b)$] $\mathbb{P}(X \geq \beta \mathbb{E}(X)) \leq \exp \{- (\beta \log \beta - \beta + 1) \mathbb{E}(X)\}$.
\end{enumerate}
\end{theorem}

Let $n$ be an odd integer, and let $G$ be an $n$-vertex Hamiltonian graph. A recipe for proving $\mathcal{C}_n(G) = \mathcal{C}(G)$ is presented in~\cite{CNP}. In order to describe it we need some definitions and results. 

The first such result stipulates the existence of a \emph{witness} subgraph of $G$; its main trait is that the incidence vector of its edge-set is orthogonal to $\mathcal{C}_n(G)$.
\begin{lemma} [\cite{CNP}] \label{lem::subgraphR}
Let $G$ be an $n$-vertex Hamiltonian graph, where $n$ is odd, and suppose that $\mathcal{C}_n(G) \neq \mathcal{C}(G)$. Then, there
exists a subgraph $R$ of $G$ such that the following conditions hold.
\begin{enumerate}
\item [\emph{(C1)}] $R \neq G$;

\item [\emph{(C2)}] Every Hamilton cycle in $G$ contains an even number of edges from $R$;

\item [\emph{(C3)}] For every partition $V (G) = A \cup B$ it holds that $e_R(A, B) \geq e_G(A, B)/2$ and $R \neq G[A, B]$.
\end{enumerate}
\end{lemma}

The following definition of a so-called \emph{parity switcher} is central to the method of~\cite{CNP}. It describes a construction that, given graphs $G$ and $R$ as in Lemma~\ref{lem::subgraphR}, aids one in finding a Hamilton cycle of $G$ with an odd number of edges in $R$, thus arriving at a contradiction to (C2) above. 

\begin{definition} \label{def::paritySwitcher}
Given a graph $G$ and a subgraph $R \subseteq G$, a subgraph $W \subseteq G$ is called an $R$-\emph{parity-switcher} if it consists of an even cycle $C = (v_1, v_2, \ldots, v_{2k}, v_1)$ with an odd number of edges in $R$, and vertex-disjoint paths $P_2, \ldots, P_k$ such that, for every $2 \leq i \leq k$, the endpoints of $P_i$ are $v_i$ and $v_{2k-i+2}$, and all internal vertices of $P_i$ are in $V(G) \setminus V(C)$.
\end{definition}

The idea behind the parity switcher gadget is that it admits two spanning paths with endpoints $v_1$ and $v_{k+1}$. Both use every edge of the connecting paths $P_i$, while one uses one perfect matching of $C$ and the other uses its other perfect matching. Therefore, the numbers of edges of $R$ on the two spanning paths have different parity.

\bigskip

We may now specify the recipe from~\cite{CNP}; it consists of the following five steps.
\begin{enumerate}

\item [(S1)] Let $G$ be an $n$-vertex Hamiltonian graph, where $n$ is odd. Suppose it satisfies $\mathcal{C}_n(G) \neq \mathcal{C}(G)$, and let $R \subseteq G$ be a witness subgraph as in Lemma~\ref{lem::subgraphR}.

\item [(S2)] Find in $G$ a, preferably small, $R$-parity-switcher $W$, that is,
\begin{enumerate}
\item [(S2a)] Find an even, preferably short, cycle $C = (v_1, \ldots, v_{2k}, v_1)$ with an odd number of edges in $R$.

\item [(S2b)] Find pairwise vertex-disjoint, preferably short, paths $P_i$ between $v_i$ and $v_{2k-i+2}$ for every $2 \leq i \leq k$.
\end{enumerate}

\item [(S3)] Find in $(G \setminus V(W)) \cup \{v_1, v_{k+1}\}$ a Hamilton path $P$ whose endpoints are $v_1$ and $v_{k+1}$.

\item [(S4)] If $P$ contains an odd (even) number of edges of $R$, then choose a Hamilton path $P'$ in $W$ whose endpoints are $v_1$ and $v_{k+1}$  with an even (odd) number of edges of $R$.

\item [(S5)] Conclude that the concatenation of $P$ and $P'$ yields a Hamilton cycle $H \subseteq G$ with an odd number of edges in R, contradicting (C2).

\end{enumerate}

Note that there is nothing to prove in steps (S4) and (S5). Moreover, whenever we start with a graph which we know to be Hamiltonian, step (S1) becomes immediate. The main task is thus to deal with steps (S2) and (S3). 

\bigskip

The following known result (see Theorem~\ref{th::HamConCexpander} below) is an important tool for handling Step (S3). In order to state it, we require the notion of an \emph{expander} and the notion of \emph{Hamilton-connectivity}.

\begin{definition} \label{def::Cexpander}
An $n$-vertex graph $G$, where $n \geq 3$, is called a $c$-expander if it satisfies the following two properties.
\begin{enumerate}
\item [\emph{(E1)}] $|N_G(X)| \geq c |X|$ holds for every $X \subseteq V(G)$ of size $|X| < n/(2c)$;

\item [\emph{(E2)}] There is an edge of $G$ between any two disjoint sets $X, Y \subseteq V(G)$ of size $|X|, |Y | \geq n/(2c)$.
\end{enumerate}
\end{definition}

\begin{definition} \label{def::HanCon}
A graph $G$ is said to be Hamilton-connected if for every two vertices $x, y \in V(G)$ there is a Hamilton path of $G$ whose endpoints are $x$ and $y$.
\end{definition}

\begin{theorem} [Theorem 7.1 in~\cite{DMMPS}] \label{th::HamConCexpander}
For every sufficiently large $c > 0$, every $c$-expander is Hamilton-connected.
\end{theorem}

The following result (see Theorem~\ref{th::PathsInExpanders} below) is our main tool for handling Step (S2b) (and is also helpful for (S3)). Before we can state it, we need the following definition of a robust expander, that is, a graph in which every small set expands even after some edges are deleted.
\begin{definition} \label{def::Palphand}
A graph $G = (V, E)$ is said to have property $P_{\alpha}(n_0, d)$ if for every $X \subseteq V$ of size $|X| \leq n_0$ and every $F \subseteq E$ such that $\left|F \cap \partial_G(x) \right| \leq \alpha \cdot \emph{deg}_G(x)$ holds for every $x \in X$, we have $|N_{G \setminus F}(X)| \geq 2 d |X|$.
\end{definition}

\begin{theorem} [Theorem 3.5 in~\cite{DKN}, abridged] \label{th::PathsInExpanders}
Let $G$ be a graph which satisfies the property $P_{\alpha}(n_0, d)$ for some $3 \leq d < n_0$. Suppose further that $e_G(A,B) > 0$ holds for any two disjoint sets $A, B \subseteq V(G)$ of sizes $|A|, |B| \geq n_0 (d-1)/16$. Let $S \subseteq V(G)$ be a set for which $|N_G(x) \cap S| \leq \beta \cdot \emph{deg}_G(x)$ holds for any vertex $x \in V(G)$. Let $a_1, \ldots, a_t, b_1, \ldots, b_t$ be $2t$ vertices in $S$, where $t \leq \frac{d n_0 \log d}{15 \log n_0}$. If $\beta < 2 \alpha - 1$, then $G$ admits pairwise vertex-disjoint paths $P_1, \ldots, P_t$ such that for every $1 \leq i \leq t$, the endpoints of $P_i$ are $a_i$ and $b_i$.
\end{theorem}

\begin{remark}
Theorem~\ref{th::PathsInExpanders} is in fact a rephrased and slightly weakened version of Theorem 3.5 from~\cite{DKN}. This version better suits our needs, as the latter theorem has an algorithmic component (as well as some other aspects) which we do not need.
\end{remark}

\section{Random graphs} \label{sec::random}

The main aim of this section is to prove Theorem~\ref{th::CycleSpaceGnp}. Before doing so, we state and prove several auxiliary results that will facilitate our proof; some are well-known and some are new. Note that the property $\mathcal{C}_n(G) = \mathcal{C}(G)$ is not monotone increasing and thus we cannot simply assume that $p$ is only slightly larger than $\log n/n$. Nevertheless, for the clarity of presentation, and since such values of $p$ are the hardest to handle, we do assume that
$$
(\log n + 2 \log \log n + \omega(1))/n \leq p := p(n) \leq 2 \log n/n .
$$
At the end of this section we briefly indicate how to adjust the proof to cope with larger values of $p$. Throughout this section it will be crucial to distinguish between vertices of typical degree and those whose degree is abnormally low; to this end we introduce the following notation: 
$$
\textrm{SMALL} := \textrm{SMALL}(G) :=  \{v \in V(G) : \textrm{deg}_G(v) \leq \log n/10\}.
$$

\bigskip

As already noted in the introduction, the smallest value of $p$ to ensure the a.a.s. Hamiltonicity of $\mathbb{G}(n,p)$ is known very precisely.
\begin{theorem} [see, e.g., Theorem 6.5 in~\cite{FK}] \label{th::HamGnp}
Let $G \sim \mathbb{G}(n,p)$, where $p := p(n) \geq (\log n + \log \log n + \omega(1))/n$. Then $G$ is a.a.s. Hamiltonian.
\end{theorem}

The threshold probability for $\mathbb{G}(n,p)$ having minimum degree at least 3 is known as well. 
\begin{theorem} [see, e.g.,~\cite{FK}] \label{th::MinDeg3}
Let $G \sim \mathbb{G}(n,p)$,  where $p := p(n) \geq (\log n + 2 \log \log n + \omega(1))/n$. Then a.a.s. $\delta(G) \geq 3$.
\end{theorem}

The following lemma lists several standard properties of random graphs. 

\begin{lemma} \label{lem::edgeDistribution}
Let $G \sim \mathbb{G}(n,p)$, where $(\log n + 2 \log \log n + \omega(1))/n \leq p := p(n) \leq 2 \log n/n$. Then, a.a.s. $G$ satisfies all of the following properties. 
\begin{enumerate}
\item [\emph{(P1)}] $\Delta(G) \leq 10 \log n$;

\item [\emph{(P2)}] $|\emph{SMALL} \cup N_G(\emph{SMALL})| \leq \sqrt{n}$;

\item [\emph{(P3)}] There is no path of length at least 1 and at most $\frac{0.3 \log n}{\log \log n}$ whose (possibly identical) endpoints lie in \emph{SMALL}.

\item [\emph{(P4)}] $e_G(A) \leq \frac{|A| \log n}{\log \log n}$ holds for every $A \subseteq [n]$ of size $|A| \leq \frac{n (\log \log n)^2}{\log n}$;

\item [\emph{(P5)}] $e_G(A, B) \leq \frac{|A| \log n}{\log \log n}$ holds for any two disjoint subsets $A, B$ of $[n]$ of size $|A| \leq \frac{n (\log \log n)^2}{\log n}$ and $|B| = |A| \sqrt{\log n}$;

\item [\emph{(P6)}] $0.999 |A| |B| p \leq e_G(A,B) \leq 1.001 |A| |B| p$ holds for any two disjoint subsets $A, B$ of $[n]$, each of size at least $\frac{n (\log \log n)^{3/2}}{\log n}$.  
\end{enumerate}
\end{lemma}

The proof of Lemma~\ref{lem::edgeDistribution} is standard and can be found in various sources (though, possibly, with slightly different parameters); for the sake of completeness we include a simple proof.
\begin{proof}
We prove that each property holds a.a.s.; since there are finitely many of them, it follows that a.a.s. they all hold simultaneously. 
\begin{enumerate}
\item [(P1)] Fix an arbitrary vertex $u \in [n]$ and note that $\textrm{deg}_G(u) \sim \textrm{Bin}(n-1, p)$. Hence
$$
\mathbb{P}(\textrm{deg}_G(u) \geq 10 \log n) \leq \binom{n-1}{10 \log n} p^{10 \log n} \leq \left(\frac{e n p}{10 \log n} \right)^{10 \log n} \leq (2 e/10)^{10 \log n} = o(1/n).
$$
A union bound over all vertices $u \in [n]$ then implies that a.a.s. $\Delta(G) \leq 10 \log n$.

\item [(P2)] Fix an arbitrary vertex $u \in [n]$ and note that $\textrm{deg}_G(u) \sim \textrm{Bin}(n-1, p)$. It thus follows by Chernoff's inequality (Theorem~\ref{th::Chernoff}(a)) that 
\begin{align*}
\mathbb{P}(u \in \textrm{SMALL}) &= \mathbb{P}(\textrm{deg}_G(u) \leq \log n/10) \leq \mathbb{P}(\textrm{deg}_G(u) \leq \mathbb{E}(\textrm{deg}_G(u))/10) \\
&\leq \exp \{- (0.1 \log 0.1 - 0.1 + 1) \log n \} \leq n^{- 0.6}.
\end{align*}
It then follows by Markov's inequality that a.a.s. $|\textrm{SMALL}| \leq n^{0.45}$. Finally, by the definition of $\textrm{SMALL}$, a.a.s. $|\textrm{SMALL} \cup N_G(\textrm{SMALL})| = O(|\textrm{SMALL}| \log n) \leq \sqrt{n}$.


\item [(P3)] Let $L = \frac{0.3 \log n}{\log \log n}$. Let $1 \leq \ell \leq L$ and let $P = (v_0, \ldots, v_{\ell})$ be a sequence of $\ell + 1$ pairwise distinct vertices, with $v_0 = v_{\ell}$ being the sole possible exception. 

Assume first that $v_0 \neq v_{\ell}$. Let $S = V(G) \setminus \{v_0, v_1, v_{\ell-1}, v_{\ell}\}$ and let ${\mathcal E}_d$ denote the event that $e_G(\{v_0, v_{\ell}\}, S) \leq \log n/5$. Note that $e_G(\{v_0, v_{\ell}\}, S) \sim \textrm{Bin}(2 |S|, p)$ and, in particular, $\mathbb{E}(e_G(\{v_0, v_{\ell}\}, S)) = 2 |S| p \geq 2 \log n$. It thus follows by Chernoff's inequality (Theorem~\ref{th::Chernoff}(a)) that
\begin{align} \label{eq::Ed}
\mathbb{P}({\mathcal E}_d) &\leq \mathbb{P}(e_G(\{v_0, v_{\ell}\}, S) \leq \mathbb{E}(e_G(\{v_0, v_{\ell}\}, S))/10) \nonumber \\
&\leq \exp \{- (0.1 \log 0.1 - 0.1 + 1) 2 \log n \} \leq n^{- 1.31}.
\end{align}
For $P$ as above, let ${\mathcal E}_P$ denote the event that $v_i v_{i+1} \in E(G)$ for every $0 \leq i < \ell$. Note that the events ${\mathcal E}_d$ and ${\mathcal E}_P$ are independent and thus $\mathbb{P}({\mathcal E}_P \wedge {\mathcal E}_d) \leq p^{\ell} n^{- 1.31}$, where the inequality is supported  by~\eqref{eq::Ed}. A union bound over all possible choices of $1 \leq \ell \leq L$ and all choices of the sequence $P = (v_0, \ldots, v_{\ell})$ implies that the probability that $G$ contains such a non-trivial path of length at most $L$ whose endpoints lie in \textrm{SMALL} is at most 
$$
\sum_{\ell=1}^{L} n^{\ell+1} p^{\ell} n^{- 1.31} \leq L (2 \log n)^L n^{- 0.31} \leq L \exp \{1.01 L \log \log n \} n^{- 0.31} \leq L n^{- 0.005} = o(1).
$$

The case $v_0 = v_{\ell}$ is similar. Let $S = V(G) \setminus \{v_1, v_{\ell-1}\}$ and let ${\mathcal E}_d$ denote the event that $\textrm{deg}_G(v_0, S) \leq \log n/10$. Similarly to the previous case, an application of Chernoff's inequality (Theorem~\ref{th::Chernoff}(a)) shows that $\mathbb{P}({\mathcal E}_d) \leq n^{- 0.6}$. Define ${\mathcal E}_P$ analogously to the previous case. Here too ${\mathcal E}_d$ and ${\mathcal E}_P$ are independent and thus $\mathbb{P}({\mathcal E}_P \wedge {\mathcal E}_d) \leq p^{\ell} n^{- 0.6}$. A union bound over all possible choices of $3 \leq \ell \leq L$ and all choices of the sequence $P = (v_0, \ldots, v_{\ell} = v_0)$ implies that the probability that $G$ contains a cycle of length at most $L$ that intersects \textrm{SMALL} is at most 
$$
\sum_{\ell=3}^L n^{\ell} p^{\ell} n^{- 0.6} = o(1).
$$
A union bound over these two cases completes the proof.

\item [(P4)] Given any set $A \subseteq [n]$ of size $1 \leq a \leq \frac{n (\log \log n)^2}{\log n}$, note that $e_G(A) \sim \textrm{Bin} \left(\binom{a}{2}, p \right)$. Hence, the probability that there exists such a set $A$ for which $e_G(A) \geq \frac{a \log n}{\log \log n}$ holds, is at most
\begin{align*}
&\sum_{a=1}^{\frac{n (\log \log n)^2}{\log n}} \binom{n}{a} \mathbb{P} \left(\textrm{Bin} \left(\binom{a}{2}, p \right) \geq \frac{a \log n}{\log \log n} \right) \\
&\leq \sum_{a=1}^{\frac{n (\log \log n)^2}{\log n}} \left(\frac{e n}{a} \left(\frac{e \binom{a}{2} p}{a \log n/\log \log n} \right)^{\log n/\log \log n} \right)^a \\
&\leq \sum_{a=1}^{\frac{n (\log \log n)^2}{\log n}} \left(\frac{e n}{a} \left(\frac{e a \log \log n}{n} \right)^{\log n/\log \log n} \right)^a \\
&\leq \sum_{a=1}^{\frac{n (\log \log n)^2}{\log n}} \left(\exp \left\{1 + \log(n/a) - \frac{\log n}{\log \log n} \left(\log(n/a) - 2 \log \log \log n \right) \right\} \right)^a = o(1).
\end{align*}

\item [(P5)] Given any set $A \subseteq [n]$ of size $1 \leq a \leq \frac{n (\log \log n)^2}{\log n}$, and any set $B \subseteq [n] \setminus A$ of size $b := a \sqrt{\log n}$, note that $e_G(A, B) \sim \textrm{Bin} \left(a^2 \sqrt{\log n}, p \right)$. Hence, the probability that there exist such sets $A$ and $B$ for which $e_G(A, B) \geq \frac{a \log n}{\log \log n}$ holds, is at most
\begin{align*}
&\sum_{a=1}^{\frac{n (\log \log n)^2}{\log n}} \binom{n}{a} \binom{n}{b} \mathbb{P} \left(\textrm{Bin} \left(a^2 \sqrt{\log n}, p \right) \geq \frac{a \log n}{\log \log n} \right) \\ 
&\leq \sum_{a=1}^{\frac{n (\log \log n)^2}{\log n}} \left(\frac{e n}{a} \left(\frac{e n}{b} \right)^{\sqrt{\log n}} \left(\frac{e a^2 \sqrt{\log n} p}{a \log n/\log \log n} \right)^{\log n/\log \log n} \right)^a \\
&\leq \sum_{a=1}^{\frac{n (\log \log n)^2}{\log n}} \left(\frac{e n}{a} \left(\frac{e n}{b} \right)^{\sqrt{\log n}} \left(\frac{6 a \sqrt{\log n} \log \log n}{n} \right)^{\log n/\log \log n} \right)^a \\
&\leq \sum_{a=1}^{\frac{n (\log \log n)^2}{\log n}} \left(\exp \left\{1 + \log (n/a) + \sqrt{\log n} (1 + \log (n/b)) - \frac{\log n}{\log \log n} (\log (n/a) - 0.6 \log \log n) \right\} \right)^a \\
&= o(1).
\end{align*}

\item [(P6)] Given a set $A \subseteq [n]$ of size $a \geq \frac{n (\log \log n)^{3/2}}{\log n}$, and a set $B \subseteq [n] \setminus A$ of size $b \geq \frac{n (\log \log n)^{3/2}}{\log n}$, note that $e_G(A, B) \sim \textrm{Bin}(a b, p)$. It thus follows by Chernoff's inequality (Theorem~\ref{th::Chernoff}) that
\begin{align*}
\mathbb{P}(e_G(A,B) \leq 0.999 a b p \textrm{ or } e_G(A,B) \geq 1.001 a b p) \leq 2 \exp \{- c a b p\} \leq \exp \{- c' a b p\},
\end{align*}
where $c, c' > 0$ are some absolute constants. A union bound over all choices of $A$ and $B$ then implies that the probability that there exist such sets $A$ and $B$ for which $0.999 |A| |B| p \leq e_G(A,B) \leq 1.001 |A| |B| p$ does not hold is at most
\begin{align*}
&\sum_{a = \frac{n (\log \log n)^{3/2}}{\log n}}^n \;\; \sum_{b = \frac{n (\log \log n)^{3/2}}{\log n}}^n \binom{n}{a} \binom{n}{b} \exp \{- c' a b p\} \\
&\leq \sum_{a = \frac{n (\log \log n)^{3/2}}{\log n}}^n \;\; \sum_{b = \frac{n (\log \log n)^{3/2}}{\log n}}^n \left(\frac{e n}{a} \right)^a \left(\frac{e n}{b} \right)^b \exp \{- c' a b p\} \\
&\leq \sum_{a = \frac{n (\log \log n)^{3/2}}{\log n}}^n \;\; \sum_{b = \frac{n (\log \log n)^{3/2}}{\log n}}^n \exp \left\{a + b + a \log (n/a) + b \log (n/b) - c' a b \log n/n \right\} \\
&\leq \sum_{a = \frac{n (\log \log n)^{3/2}}{\log n}}^n \;\; \sum_{b = \frac{n (\log \log n)^{3/2}}{\log n}}^n \exp \left\{a \log \log n + b \log \log n - c' a (\log \log n)^{3/2}/2 - c' b (\log \log n)^{3/2}/2 \right\} \\
&= o(1).
\end{align*}
\end{enumerate}
\end{proof}

The following lemma asserts that for any two disjoint vertex sets of appropriate sizes in $\mathbb{G}(n,p)$, some vertex in one set must have a large degree in the other.
\begin{lemma} \label{lem::halfDegree}
Let $\delta \in (0, 1/4)$ be an arbitrarily small constant, and let $G \sim \mathbb{G}(n,p)$, where $p := p(n) \geq (\log n + 2 \log \log n + \omega(1))/n$. Then a.a.s. the following holds. For every $A \subseteq V(G)$ of size $|A| \geq \frac{n (\log \log n)^2}{\sqrt{\log n}}$ and every $B \subseteq [n] \setminus A$ of size $|B| = (1/2 + \delta) n$ there exists a vertex $u \in A$ such that $\emph{deg}_G(u, B) \geq (1 + \delta) \emph{deg}_G(u)/2$.
\end{lemma}

\begin{proof}
Observe that it suffices to prove the lemma for all sets $A$ whose size is precisely $a := \frac{n (\log \log n)^2}{\sqrt{\log n}}$. A standard application of Chernoff's inequality (Theorem~\ref{th::Chernoff}) and a union bound shows that a.a.s. $G$ satisfies all of the following properties.
\begin{enumerate}
\item [(Q1)] $e_G(A) \leq a^2 p$ holds for every $A \subseteq [n]$ of size $a$.

\item [(Q2)] $(1 - \delta/4) a (1/2 + \delta) n p \leq e_G(A,B) \leq (1 + \delta/4) a (1/2 + \delta) n p$ holds for every $A \subseteq [n]$ of size $a$ and every $B \subseteq [n] \setminus A$ of size $|B| = (1/2 + \delta) n$.

\item [(Q3)] $e_G(A,[n] \setminus (A \cup B)) \leq (1 + \delta/4) a (1/2 - \delta) n p$ holds for every $A \subseteq [n]$ of size $a$ and every $B \subseteq [n] \setminus A$ of size $|B| = (1/2 + \delta) n$.
\end{enumerate} 
Assume then that $G$ satisfies properties (Q1), (Q2), and (Q3). Suppose for a contradiction that there exist sets $A \subseteq [n]$ of size $a$ and $B \subseteq [n] \setminus A$ of size $|B| = (1/2 + \delta) n$ such that $\textrm{deg}_G(u, B) < (1 + \delta) \textrm{deg}_G(u)/2$ holds for every $u \in A$. It then follows that
\begin{align*}
(1 - \delta/4) a (1/2 + \delta) n p \leq e_G(A, B) = \sum_{u \in A} \textrm{deg}_G(u, B) < \sum_{u \in A} (1 + \delta) \textrm{deg}_G(u)/2,
\end{align*}
implying that $\sum_{u \in A} \textrm{deg}_G(u) > (1 + \delta/3) a n p$. On the other hand
\begin{align*}
\sum_{u \in A} \textrm{deg}_G(u) &= \sum_{u \in A} \textrm{deg}_G(u, A) + \sum_{u \in A} \textrm{deg}_G(u, B) + \sum_{u \in A} \textrm{deg}_G(u, [n] \setminus (A \cup B)) \\
&= 2 e_G(A) + e_G(A,B) + e_G(A, [n] \setminus (A \cup B)) \leq (1 + \delta/3) a n p,
\end{align*}
which is a clear contradiction.
\end{proof}

The following result allows us to split a graph into several parts in a beneficial manner. The specific formulation we use is taken from~\cite{HKT}, though similar results can be found in other sources.
\begin{lemma} [Lemma 2.4 in~\cite{HKT}] \label{lem::LLLsplit}
Let $G = (V, E)$ be a graph on $n$ vertices with maximum degree $\Delta$. Let $Y \subseteq V$ be a set of $m = a + b$ vertices, where $a$ and $b$ are positive integers. Assume that $\emph{deg}_G(v, Y) \geq \delta$ holds for every $v \in V$. If $\Delta^2 \cdot \lceil \frac{m}{\min \{a,b\}} \rceil \cdot 2 \cdot e^{1 - \frac{\min \{a,b\}^2}{5 m^2} \delta} < 1$, then there exists a partition $Y = A \cup B$ of $Y$ such that the following properties hold.
\begin{enumerate}
\item [$(1)$] $|A| = a$ and $|B| = b$;

\item [$(2)$] $\emph{deg}_G(v, A) \geq \frac{a}{3m} \emph{deg}_G(v, Y)$ holds for every $v \in V$;

\item [$(3)$] $\emph{deg}_G(v, B) \geq \frac{b}{3m} \emph{deg}_G(v, Y)$ holds for every $v \in V$.  
\end{enumerate}
\end{lemma}

The essential implication of the following result is that high degree subgraphs of random graphs are good expanders.

\begin{lemma} \label{lem::subgraphExpander}
Let $G$ be an $n$-vertex graph which satisfies properties \emph{(P4)} and \emph{(P5)} from Lemma~\ref{lem::edgeDistribution}. Let $H$ be a (not necessarily spanning) subgraph of $G$ for which $\delta(H) \geq \gamma \log n$ holds for some constant $\gamma > 0$. Let $\alpha \in [0,1)$ be a constant. Then, $|N_{H \setminus F}(X)| \geq |X| \sqrt{\log n}$ for every $F \subseteq E(H)$ and every $X \subseteq V(H)$ of size $|X| \leq \frac{n (\log \log n)^2}{\log n}$ such that $\left|F \cap \partial_H(x) \right| \leq \alpha \cdot \emph{deg}_H(x)$ holds for every $x \in X$. 
\end{lemma}

\begin{proof}
Suppose for a contradiction that $X \subseteq V(H)$ is a set of size $|X| \leq \frac{n (\log \log n)^2}{\log n}$ and $F \subseteq E(H)$ is such that $\left|F \cap \partial_H(x) \right| \leq \alpha \cdot \textrm{deg}_H(x)$ holds for every $x \in X$, and yet $|N_{H \setminus F}(X)| < |X| \sqrt{\log n}$. It follows by the premise of the lemma that
$$
|X| \cdot (1 - \alpha) \gamma \log n \leq \sum_{x \in X} \textrm{deg}_{H \setminus F}(x) \leq 2 e_G(X) + e_G(X, N_{H \setminus F}(X)) \leq \frac{3 |X| \log n}{\log \log n}, 
$$
which is an obvious contradiction. 
\end{proof}

The following lemma asserts that if a subgraph of $\mathbb{G}(n,p)$ inherits at least half of its edges in every cut, then it also inherits at least one edge between any two large disjoint sets. 
\begin{lemma} \label{lem::robustEdgeBetweenSets}
Let $G \sim \mathbb{G}(n,p)$, where $(\log n + 2 \log \log n + \omega(1))/n \leq p := p(n) \leq 2 \log n/n$. Then a.a.s. the following holds. Let $R$ be a subgraph of $G$ such that $e_R(A, V(G) \setminus A) \geq e_G(A, V(G) \setminus A)/2$ for every $A \subseteq V(G)$. Then, $e_R(A, B) > 0$ holds for any two disjoint sets $A, B \subseteq V(G)$ of size $|A| = |B| = 2n/5$. 
\end{lemma}

\begin{proof}
Suppose that $G$ satisfies Property (P6) from Lemma~\ref{lem::edgeDistribution} (this fails with probability $o(1)$). Fix any two disjoint sets $A, B \subseteq V(G)$ of size $|A| = |B| = 2n/5$, and let $D = V(G) \setminus (A \cup B)$. It follows by Property (P6) that $e_G(A, V(G) \setminus A) \geq 0.999 |A| |V(G) \setminus A| p \geq 0.999 \cdot \frac{6}{25} \cdot n^2 p$ and that $e_G(A, D) \leq 1.001 |A| |D| p \leq 1.001 \cdot \frac{2}{25} \cdot n^2 p$. Therefore,
\begin{align*}
e_R(A, B) &= e_R(A, V(G) \setminus A) - e_R(A,D) \geq e_G(A, V(G) \setminus A)/2 - e_G(A,D) \\
&\geq 0.999 \cdot \frac{3}{25} \cdot n^2 p - 1.001 \cdot \frac{2}{25} \cdot n^2 p > 0.
\end{align*}
\end{proof}

When building a parity switcher in Step (S2), it is useful to keep it small. The following result is helpful in that respect.

\begin{lemma} \label{lem::DiameterR}
Let $G \sim \mathbb{G}(n,p)$, where $(\log n + 2 \log \log n + \omega(1))/n \leq p := p(n) \leq 2 \log n/n$. Then a.a.s. the following holds. Let $R$ be a subgraph of $G$ such that $e_R(A, V(G) \setminus A) \geq e_G(A, V(G) \setminus A)/2$ holds for every $A \subseteq V(G)$. Then, for every set $S \subseteq V(G)$ of size $|S| = o(\log n)$ and every two vertices $x, y \in V(G) \setminus (\emph{SMALL} \cup S)$, there is a path between $x$ and $y$ in $(R \setminus (\emph{SMALL} \cup N_G(\emph{SMALL}) \cup S)) \cup \{x,y\}$ whose length is at most $5 \frac{\log n}{\log \log n}$.
\end{lemma}

\begin{proof}
Fix an arbitrary set $S \subseteq V(G)$ of size $|S| = o(\log n)$, and let $Z = \textrm{SMALL} \cup N_G(\textrm{SMALL}) \cup S$; note that $|Z| \leq 2 \sqrt{n}$ holds a.a.s. by Property (P2) from Lemma~\ref{lem::edgeDistribution}. Note that a.a.s. $\textrm{deg}_G(u, \textrm{SMALL} \cup N_G(\textrm{SMALL})) \leq 1$ holds for every $u \in V(G) \setminus \textrm{SMALL}$ by Property (P3) from Lemma~\ref{lem::edgeDistribution}. Since, moreover, $|S| = o(\log n)$ and $\deg_R(v) \geq \deg_G(v)/2$ holds for every $v \in V(G)$ by the premise of the lemma, it follows that $\textrm{deg}_{R \setminus Z}(x) \geq \log n/21$ holds for every $x \in V(G) \setminus Z$. For every vertex $x \in V(G) \setminus Z$ and every non-negative integer $i$, let $N_{R \setminus Z}^i(x) = \{u \in V(G) \setminus Z : \textrm{dist}_{R \setminus Z}(x,u) \leq i\}$. Starting from an arbitrary vertex $x \in V(G) \setminus Z$, repeated applications of Lemma~\ref{lem::subgraphExpander} with $\alpha = 0$ and $H = R \setminus Z$ show that a.a.s. there exists an integer $t \leq (2 + o(1)) \frac{\log n}{\log \log n}$ such that $\left|N_{R \setminus Z}^t(x) \right| \geq \frac{n (\log \log n)^2}{\sqrt{\log n}}$. We claim that a.a.s. $\left|N_{R \setminus Z}^{t+1}(x) \right| \geq 0.48 n$. Indeed, if not, then there exists a set $B \subseteq V(G) \setminus \left(Z \cup N_{R \setminus Z}^{t+1}(x) \right)$ of size $|B| = 0.51 n$ such $E_{R \setminus Z} \left(N_{R \setminus Z}^t(x), B \right) = \varnothing$. Observe that a.a.s. $\textrm{deg}_{R \setminus Z}(v) \geq (1/2 - o(1)) \textrm{deg}_G(v)$ holds for every vertex $v \in V(G) \setminus Z$. Indeed, $v \notin \textrm{SMALL}$, $\textrm{deg}_R(v) \geq \textrm{deg}_G(v)/2$, $|S| = o(\log n)$, and $\textrm{deg}_R(v, \textrm{SMALL} \cup N_G(\textrm{SMALL})) \leq 1$ holds a.a.s. by Property (P3) from Lemma~\ref{lem::edgeDistribution}. It thus follows that a.a.s. $\textrm{deg}_G(v, B) < (1/2 + o(1)) \textrm{deg}_G(v)$ holds for every $v \in N_{R \setminus Z}^t(x)$. However, by Lemma~\ref{lem::halfDegree}, this occurs with probability $o(1)$. 


Let $x, y \in V(G) \setminus (\textrm{SMALL} \cup S)$ be any two vertices. If $x \in N_G(\textrm{SMALL})$, then let $x' \in N_R(x) \setminus Z$ be arbitrary; otherwise, set $x' = x$. Similarly, if $y \in N_G(\textrm{SMALL})$, then let $y' \in N_R(y) \setminus Z$ be arbitrary; otherwise, set $y' = y$. Since $x', y' \in V(G) \setminus Z$, the above argument implies that a.a.s. $\left|N_{R \setminus Z}^{t+1}(x') \right| \geq 0.48 n$ and $\left|N_{R \setminus Z}^{t+1}(y') \right| \geq 0.48 n$. If $N_{R \setminus Z}^{t+1}(x') \cap N_{R \setminus Z}^{t+1}(y') \neq \varnothing$, then there is a path of length at most $2t+2$ between $x'$ and $y'$ in $R \setminus Z$. Otherwise, it follows by Lemma~\ref{lem::robustEdgeBetweenSets} that a.a.s. there is an edge of $R$ between $N_{R \setminus Z}^{t+1}(x')$ and $N_{R \setminus Z}^{t+1}(y')$, yielding a path of length at most $2t+3$ between $x'$ and $y'$ in $R \setminus Z$. In either case, adding the edges $x x'$ and $y y'$, if they exist, yields a path of length at most $2t+5 \leq 5 \frac{\log n}{\log \log n}$ between $x$ and $y$ in $(R \setminus Z) \cup \{x,y\}$. 
\end{proof}

Lemma 2.3 in~\cite{CNP} is the main tool for handling Step (S2a) in that paper. Unfortunately, using this lemma as a black box might create a cycle which exhausts the neighbourhood of some vertex outside the cycle; that vertex cannot be absorbed in any of the following steps. In order to circumvent this problem, we state and prove a variant of Lemma 2.3 from~\cite{CNP} that is suitable for random graphs having vertices of very small degrees.
\begin{lemma} \label{lem::OurCycleLemma}
Let $G \sim \mathbb{G}(n,p)$, where $(\log n + 2 \log \log n + \omega(1))/n \leq p := p(n) \leq 2 \log n/n$. Then a.a.s. the following holds. Let $R \neq G$ be a subgraph of $G$ such that $e_R(A, V(G) \setminus A) \geq e_G(A, V(G) \setminus A)/2$ for every $A \subseteq V(G)$, and $R \neq G[A, B]$ for every partition $V(G) = A \cup B$. Then, there exists a cycle $C \subseteq G$ satisfying all of the following properties.
\begin{enumerate}
\item [$(a)$] $|C|$ is even and $|E(C) \setminus E(R)| = 1$;

\item [$(b)$] $|C| \leq 22 \frac{\log n}{\log \log n}$;

\item [$(c)$] $\emph{deg}_G(u, V(C)) \leq 2$ holds for every $u \in \emph{SMALL} \cap V(C)$;

\item [$(d)$] $\emph{deg}_G(u, V(C)) \leq 1$ holds for every $u \in \emph{SMALL} \setminus V(C)$.
\end{enumerate}
\end{lemma}

\begin{proof}
Suppose that $G$ satisfies the assertions of Theorem~\ref{th::MinDeg3}, and of Lemmas~\ref{lem::edgeDistribution} and~\ref{lem::DiameterR}; note that this fails with probability $o(1)$. 

Let $Z = \textrm{SMALL} \cup N_G(\textrm{SMALL}) $ and let $U = V(G) \setminus Z$. It follows by Lemma~\ref{lem::DiameterR} that $R[U]$ is connected. We distinguish between the following two cases.
\begin{enumerate}
\item [(1)] $R[U]$ is bipartite. Let $U = A \cup B$ be the unique bipartition of $R[U]$. This case is divided further into the following subcases.
\begin{enumerate}
\item [(1.1)] There exist vertices $x \in A$ and $y \in B$ such that $xy \in E(G) \setminus E(R)$. It follows by Lemma~\ref{lem::DiameterR} that there is a path in $R[U]$ between $x$ and $y$ of length $\ell \leq 5 \frac{\log n}{\log \log n}$. Since $R[U]$ is bipartite, $\ell$ is odd. Hence, combined with the edge $xy$, this yields the required cycle $C$ (properties (c) and (d) are satisfied since $V(C) \cap Z = \varnothing$).

\item [(1.2)] $E_G(A,B) = E_R(A,B)$. This case is again divided into two subcases.
\begin{enumerate}
\item [(1.2.1)] $R$ is bipartite. It follows by Theorem~\ref{th::MinDeg3}, since $\textrm{deg}_R(v) \geq \textrm{deg}_G(v)/2$ holds for every $v \in V(G)$ by the premise of the lemma, by Property (P3) from Lemma~\ref{lem::edgeDistribution}, and by Lemma~\ref{lem::DiameterR} that $R$ is connected (indeed, by Lemma~\ref{lem::DiameterR}, $R \setminus \textrm{SMALL}$ is connected, and by the other aforementioned results, any $u \in \textrm{SMALL}$ satisfies $N_R(u, V(G) \setminus \textrm{SMALL}) \neq \emptyset$). Let $V(G) = A' \cup B'$ be the unique bipartition of $R$; note that, up to renaming the parts, $A \subseteq A'$ and $B \subseteq B'$ must hold. Since $R \neq G[A', B']$ by the premise of the lemma, there exist vertices $x \in A'$ and $y \in B'$ such that $xy \in E(G) \setminus E(R)$. Since $E_G(A,B) = E_R(A,B)$ it follows that $\{x,y\} \cap Z \neq \varnothing$; assume without loss of generality that $x \in Z$. Assume first that $x \in \textrm{SMALL}$ and note that $y \notin \textrm{SMALL}$ by Property (P3) from Lemma~\ref{lem::edgeDistribution}. Let $y' \in B' \setminus \{y\}$ be a vertex for which $xy' \in E(R)$ and note that $y' \notin \textrm{SMALL}$ holds by Property (P3) from Lemma~\ref{lem::edgeDistribution}. It then follows by Lemma~\ref{lem::DiameterR} that there is a path in $R[U \cup \{y, y'\}]$ between $y$ and $y'$ of length $\ell \leq 5 \frac{\log n}{\log \log n}$. Since $R$ is bipartite, $\ell$ is even. Hence, combined with the edges $xy$ and $xy'$, this yields a cycle $C$ satisfying properties (a) and (b). Moreover, properties (c) and (d) are satisfied by Property (P3) from Lemma~\ref{lem::edgeDistribution} and since $V(C) \cap \textrm{SMALL} = \{x\}$ and $V(C) \cap N_G(\textrm{SMALL}) = \{y, y'\}$. Assume then that $x \in N_G(\textrm{SMALL})$. If $y \in \textrm{SMALL}$, then this is analogous to the previous case; by Property (P3) from Lemma~\ref{lem::edgeDistribution} we may thus assume that $y \in U$. It then follows by Lemma~\ref{lem::DiameterR} that there is a path in $R[U \cup \{x\}]$ between $x$ and $y$ of length $\ell \leq 5 \frac{\log n}{\log \log n}$. Since $R$ is bipartite, $\ell$ is odd. Hence, combined with the edge $xy$, this yields the required cycle $C$ (properties (c) and (d) are satisfied since $V(C) \cap \textrm{SMALL} = \varnothing$ and $V(C) \cap N_G(\textrm{SMALL}) = \{x\}$).

\item [(1.2.2)] $R$ is not bipartite. Let $C'$ be an odd cycle in $R$ and note that $V(C') \cap Z \neq \varnothing$. Assume first that $V(C') \cap \textrm{SMALL} = \varnothing$. Since, moreover, $N_G(\textrm{SMALL})$ is an independent set by Property (P3) from Lemma~\ref{lem::edgeDistribution}, it follows that $N_{C'}(u) \subseteq U$ holds for every $u \in V(C') \cap N_G(\textrm{SMALL})$. It is then straightforward to verify that there must exist a vertex $w \in V(C') \cap N_G(\textrm{SMALL})$ such that $\textrm{deg}_{C'}(w, A) = 1$ and $\textrm{deg}_{C'}(w, B) = 1$. Let $x$ be the unique element of $N_{C'}(w, A)$ and let $y$ be the unique element of $N_{C'}(w, B)$. Note that $G[U]$ is a.a.s. not bipartite (indeed, otherwise, both $A$ and $B$ are independent in $G$, but by Property (P2) from Lemma~\ref{lem::edgeDistribution}, at least one of them is of size $\Omega(n)$; by Property (P6) no such set exists). Since $R[U]$ is bipartite, we may assume without loss of generality that there are vertices $u, v \in A$ such that $uv \in E(G) \setminus E(R)$. Since $U \cap \textrm{SMALL} = \emptyset$, it follows by Lemma~\ref{lem::DiameterR} that there is a path $P_x$ in $R[U \setminus \{y,v\}]$ between $x$ and $u$ of length $\ell_1 \leq 5 \frac{\log n}{\log \log n}$. Since $R[U]$ is bipartite, $\ell_1$ is even. Similarly, since $|V(P_x)| = o(\log n)$, it follows by Lemma~\ref{lem::DiameterR} that there is a path $P_y$ in $R[U] \setminus V(P_x)$ between $y$ and $v$ of length $\ell_2 \leq 5 \frac{\log n}{\log \log n}$. Since $R[U]$ is bipartite, $\ell_2$ is odd. Then $w x P_x u v P_y y w$ forms the required cycle $C$ (properties (c) and (d) are satisfied since $V(C) \cap Z = \{w\}$).

The case $V(C') \cap \textrm{SMALL} \neq \varnothing$ is essentially the same. The only difference is that either there exists a vertex $w \in V(C') \cap N_G(\textrm{SMALL})$ with one neighbour in $V(C') \cap A$ and the other in $V(C') \cap B$ as in the previous case, or there is a vertex $w \in V(C') \cap \textrm{SMALL}$ for which there exist two vertices $x', y' \in N_{C'}(w)$ such that $x' x \in E(R)$ for some $x \in A$ and $y' y \in E(R)$ for some $y \in B$. The obtained cycle $C$ is then $w x' x P_x u v P_y y y' w$; note that all of these edges except for $uv$ are in $R$.
\end{enumerate}
\end{enumerate}

\item [(2)] $R[U]$ is not bipartite. Let $C' = (x_1, \ldots, x_{2t-1}, x_1)$ be a shortest odd cycle in $R[U]$. Note that $|C'| \leq 11 \frac{\log n}{\log \log n}$. Indeed, suppose for a contradiction that $|C'| > 11 \frac{\log n}{\log \log n}$. By Lemma~\ref{lem::DiameterR} there is a path $P$ in $R[U]$ between $x_1$ and $x_t$ whose length is at most $5 \frac{\log n}{\log \log n}$. However $P \cup C'$ contains an odd cycle which is shorter than $C'$, contrary to the assumed minimality of $|C'|$ (note that $P$ and $C'$ may have several vertices and edges in common, but this does not invalidate the conclusion).

Since $R \neq G$ by the premise of the lemma, there exists an edge $xy \in E(G) \setminus E(R)$. By Property (P3) from Lemma~\ref{lem::edgeDistribution} we may assume without loss of generality that $y \notin \textrm{SMALL}$. We claim that there exists a vertex $u \in V(C') \setminus \{y\}$ such that there exists a path $P_x$ in $(R \setminus (V(C') \cup \{y\})) \cup \{x, u\}$ between $x$ and $u$ of length at most $5 \frac{\log n}{\log \log n} + 1$ and such that $|V(P_x) \cap \textrm{SMALL}| \leq 1$ and $|V(P_x) \cap N_G(\textrm{SMALL})| \leq 1$. If $x \in V(C')$ or $N_R(x) \cap V(C') \neq \varnothing$, then this is obvious. Otherwise, let $x' \in N_R(x) \setminus \textrm{SMALL}$ be an arbitrary vertex if $x \in \textrm{SMALL}$ and let $x' = x$ if $x \notin \textrm{SMALL}$. Note that such a vertex $x'$ exists as $\textrm{deg}_R(x) \geq \textrm{deg}_G(x)/2$ holds by the premise of the lemma, since $\delta(G) \geq 3$ holds by Theorem~\ref{th::MinDeg3}, and due to Property (P3) from Lemma~\ref{lem::edgeDistribution}. Since $|V(C') \cup \{y\}| = o(\log n)$, by Lemma~\ref{lem::DiameterR}, there is a path $P$ in $(R \setminus Z) \cup \{x'\}$ between $x'$ and some vertex $u \in V(C')$ of length at most $5 \frac{\log n}{\log \log n}$, such that $V(P) \cap (V(C') \cup \{y\}) = \{u\}$. Adding, if needed, the edge $x x'$, yields the required path $P_x$. Similarly, there exists a path $P_y$ in $(R[U] \setminus (V(P_x) \cup V(C'))) \cup \{y, v\}$ between $y$ and a vertex $v \in V(C') \setminus \{u\}$ of length at most $5 \frac{\log n}{\log \log n}$. Indeed, if $y \in V(C')$, then this is obvious. Otherwise, let $v \in V(C') \setminus \{u\}$ be an arbitrary vertex; note that $v \notin \textrm{SMALL}$. Since $y \notin \textrm{SMALL}$ and $|V(P_x) \cup V(C')| = o(\log n)$, the required path $P_y$ exists by Lemma~\ref{lem::DiameterR}. Combining the path $u P_x x y P_y v$ (that has precisely one edge in $E(G) \setminus E(R)$) with one of the two paths that connect $u$ and $v$ in $C'$ (all of whose edges are in $R$) yields a cycle $C$ satisfying properties (a) and (b). Moreover, properties (c) and (d) are satisfied since either $V(C) \cap \textrm{SMALL} = \varnothing$ and $|V(C) \cap N_G(\textrm{SMALL})| \leq 1$ or $V(C) \cap \textrm{SMALL} = \{x\}$ and $V(C) \cap N_G(\textrm{SMALL}) = \{x', y\}$.
\end{enumerate}
\end{proof} 

The next result is our main tool for handling Step (S3). 
\begin{lemma} \label{lem::NewHamCon}
Let $G \sim \mathbb{G}(n,p)$, where $(\log n + 2 \log \log n + \omega(1))/n \leq p := p(n) \leq 2 \log n/n$. Then a.a.s. the following holds. Let $S \subseteq V(G)$ be a set of size $|S| = \Omega(n)$, let $x, y \in S \setminus \emph{SMALL}$ be any two vertices, and let $G' = G[S]$. Suppose that $\emph{deg}_G(u, S) \geq \gamma \log n$ holds for some constant $\gamma > 0$ and every $u \in S \setminus \emph{SMALL}$, and that $\emph{deg}_G(u, S \setminus \{x,y\}) \geq 2$ holds for every $u \in S \cap \emph{SMALL}$. Then there exists a Hamilton path of $G'$ whose endpoints are $x$ and $y$.
\end{lemma}

\begin{proof}
Suppose that $G$ satisfies the assertion of Lemma~\ref{lem::edgeDistribution}; note that this fails with probability $o(1)$.

Let $S \cap \textrm{SMALL} = \{u_1, \ldots, u_t\}$; suppose first that $t \geq 1$. Let $x_1, y_1, \ldots, x_t, y_t$ be $2t$ distinct vertices such that $x_i, y_i \in N_{G'}(u_i, S) \setminus (\textrm{SMALL} \cup \{x,y\})$ for every $1 \leq i \leq t$; such vertices exist by the premise of the lemma and by Property (P3) from Lemma~\ref{lem::edgeDistribution}.

Let $U = \{u_1, \ldots, u_t\} \cup \{x_1, y_1, \ldots, x_t, y_t\}$ and note that $\textrm{deg}_G(u, U) \leq 1$ holds for every $u \in S \setminus \textrm{SMALL}$ by Property (P3) from Lemma~\ref{lem::edgeDistribution}. Since, moreover,  $G$ satisfies Property (P1) from Lemma~\ref{lem::edgeDistribution}, we may apply Lemma~\ref{lem::LLLsplit} to $G[S \setminus \textrm{SMALL}]$ and obtain a partition $S_1 \cup S_2$ of $Y := S \setminus (U \cup \{x,y\})$ into two parts of essentially equal size such that $\textrm{deg}_G(u, S_1) \geq \gamma \log n/10$ and $\textrm{deg}_G(u, S_2) \geq \gamma \log n/10$ hold for every $u \in S \setminus \textrm{SMALL}$. Let $G_1 = G[S_1 \cup \{y_1, x_2, \ldots, x_t, y_t, y\}]$ and note that, by Lemma~\ref{lem::subgraphExpander}, the graph $G_1$ satisfies the property $P_{2/3}(n/\log n, \sqrt{\log n}/2)$. Moreover, by Property (P6) of Lemma~\ref{lem::edgeDistribution}, there is an edge of $G_1$ between any two disjoint subsets of $V(G_1)$, each of size at least $n/(40 \sqrt{\log n})$. Setting $L = \{y_1, x_2, \ldots, x_t, y_t, y\}$, observe that $|N_{G_1}(z) \cap L| \leq 2 \leq \textrm{deg}_{G_1}(z)/10$ holds for any vertex $z \in V(G_1)$. It thus follows by Theorem~\ref{th::PathsInExpanders} that $G_1$ admits pairwise vertex-disjoint paths $P_1, P_2, \ldots, P_t$, where the endpoints of $P_1$ are $y_t$ and $y$, and for every $2 \leq i \leq t$, the endpoints of $P_i$ are $x_i$ and $y_{i-1}$. 

Let $W = S \setminus (V(P_1) \cup \ldots \cup V(P_t) \cup \{u_1, \ldots, u_t\})$ and let $G_2 = G[W \cup \{x\}]$. Since $S_2 \subseteq W$ and $W \cap \textrm{SMALL} = \varnothing$, it follows that $\delta(G_2) \geq \gamma \log n/10$. Hence, it follows by Property (P6) from Lemma~\ref{lem::edgeDistribution} and by Lemma~\ref{lem::subgraphExpander} that $G_2$ is a $c$-expander, where $c$ is a sufficiently large constant, as per Theorem~\ref{th::HamConCexpander} (note that the expansion of large sets, which are not covered by Lemma~\ref{lem::subgraphExpander}, is ensured by Property (P6)). It then follows by Theorem~\ref{th::HamConCexpander} that $G_2$ admits a Hamilton path whose endpoints are $x$ and $x_1$. Combined with the previously built paths $P_1, P_2, \ldots, P_t$ and with the paths $x_1 u_1 y_1, \ldots, x_t u_t y_t$, this yields a Hamilton path of $G'$ whose endpoints are $x$ and $y$.

Finally, if $t=0$, then $G'$ itself is a.a.s. a $c$-expander and is thus a.a.s. Hamilton-connected by Theorem~\ref{th::HamConCexpander}.
\end{proof}

We are now in a position to prove Theorem~\ref{th::CycleSpaceGnp}. 

\begin{proof} [Proof of Theorem~\ref{th::CycleSpaceGnp}]
Let $G \sim \mathbb{G}(n,p)$, where $(\log n + 2 \log \log n + \omega(1))/n \leq p := p(n) \leq 2 \log n/n$. Suppose that $G$ satisfies the assertions of Theorems~\ref{th::HamGnp} and~\ref{th::MinDeg3}, and of Lemma~\ref{lem::edgeDistribution}; note that this fails with probability $o(1)$. 

Suppose for a contradiction that $\mathcal{C}_n(G) \neq \mathcal{C}(G)$. We follow the recipe that was presented in Section~\ref{sec::prelim}. That is, we need to handle steps (S1), (S2), and (S3). Our assumption that $\mathcal{C}_n(G) \neq \mathcal{C}(G)$ will then lead to the contradiction appearing in (S5).  

 Let $R$ be a witness subgraph of $G$ as in the premise of Lemma~\ref{lem::subgraphR}; combined with Theorem~\ref{th::HamGnp}, this takes care of (S1). Suppose that $G$  satisfies the assertions of Lemmas~\ref{lem::DiameterR} and~\ref{lem::OurCycleLemma} (with respect to $R$); note that this fails with probability $o(1)$. 

Next, we take care of (S2). Starting with (S2a), it follows by Lemma~\ref{lem::OurCycleLemma} that $G$ contains a cycle $C$ satisfying properties (a), (b), (c), and (d). In particular, $C = (v_1, \ldots, v_{2k})$ is an even cycle having an odd number of edges in $R$, and its length is at most $22 \frac{\log n}{\log \log n}$.

Prior to handling (S2b) and thinking ahead to Step (S3), let $v'_1, \ldots, v'_{2k}$ be distinct vertices, where, for every $1 \leq i \leq 2k$, $v'_i \in N_G(v_i) \setminus (\textrm{SMALL} \cup V(C))$ if $v_i \in \textrm{SMALL}$ and $v'_i = v_i$ if $v_i \notin \textrm{SMALL}$; such vertices $v'_1, \ldots, v'_{2k}$ exist by Theorem~\ref{th::MinDeg3}, by Property (c) from Lemma~\ref{lem::OurCycleLemma}, and by Property (P3) from Lemma~\ref{lem::edgeDistribution}. Let $U = \{v_1, \ldots, v_{2k}\} \cup \{v'_1, \ldots, v'_{2k}\}$, and note that $|U| =  o(\log n)$. Moreover, note that similarly to the argument for the existence of $v'_1, \ldots, v'_{2k}$, it follows by Theorem~\ref{th::MinDeg3}, by Property (d) from Lemma~\ref{lem::OurCycleLemma}, and by Property (P3) from Lemma~\ref{lem::edgeDistribution} that $\textrm{deg}_G(u, V(G) \setminus U) \geq 2$ for every $u \in \textrm{SMALL} \setminus V(C)$.  

Let $Z = \textrm{SMALL} \cup N_G(\textrm{SMALL})$, and note that $\textrm{deg}_G(u, Z) \leq 1$ holds for every $u \in V(G) \setminus \textrm{SMALL}$ by Property (P3) from Lemma~\ref{lem::edgeDistribution}. Since, moreover, $|U| = o(\log n)$, by Property (P1) from Lemma~\ref{lem::edgeDistribution} we may apply Lemma~\ref{lem::LLLsplit} to $G[V(G) \setminus \textrm{SMALL}]$ with $Y = V(G) \setminus (\textrm{SMALL} \cup U)$ to obtain a set $A' \subseteq V(G) \setminus (\textrm{SMALL} \cup U)$ of size $n/2$ and a set $B' := V(G) \setminus (A' \cup \textrm{SMALL})$ such that $\textrm{deg}_G(u, A') \geq \log n/65$ and $\textrm{deg}_G(u, B') \geq \log n/65$ hold for every $u \in V(G) \setminus \textrm{SMALL}$. Let $A = (A' \cup Z) \setminus U$ and let $B = (B' \setminus Z) \cup \{v'_1, \ldots, v'_{2k}\}$; note that, by properties (P2) and (P3) from Lemma~\ref{lem::edgeDistribution}, $|A| = (1/2 \pm o(1)) n$, $|B| = (1/2 \pm o(1)) n$, and, moreover, $\textrm{deg}_G(u, A) \geq \log n/70$ and $\textrm{deg}_G(u, B) \geq \log n/70$ hold for every $u \in V(G) \setminus \textrm{SMALL}$.   

Returning to (S2b), let $G_1$ be the graph obtained from $G[B]$ by deleting $v'_1$ and $v'_{k+1}$ and all the edges (but none of the other vertices) of $E(C)$; note that $V(G_1) \cap \textrm{SMALL} = \varnothing$ and thus $\delta(G_1) \geq \log n/70 - 4 \geq \log n/71$. It thus follows by Lemma~\ref{lem::subgraphExpander} that $G_1$ satisfies the property $P_{2/3}(n/\log n, \sqrt{\log n}/2)$. Moreover, by Property (P6) of Lemma~\ref{lem::edgeDistribution}, there is an edge of $G_1$ between any two disjoint subsets of $V(G_1)$, each of size at least $n/(40 \sqrt{\log n})$. Setting $S = \{v'_2, \ldots, v'_k, v'_{k+2}, \ldots, v'_{2k}\}$ and noting that $|S| = o(\delta(G_1))$, observe that $|N_{G_1}(x) \cap S| \leq \textrm{deg}_{G_1}(x)/10$ holds for any vertex $x \in V(G_1)$. It thus follows by Theorem~\ref{th::PathsInExpanders} that $G_1$ admits pairwise vertex-disjoint paths $P'_2, \ldots, P'_k$ such that, for every $2 \leq i \leq k$, the endpoints of $P'_i$ are $v'_i$ and $v'_{2k-i+2}$. Finally, for every $2 \leq i \leq k$, let $P_i = v_i v'_i P'_i v'_{2k-i+2} v_{2k-i+2}$ (where $v_j v'_j$ is simply the vertex $v_j$ if $v'_j = v_j$); note that $P_2, \ldots, P_k$ are pairwise vertex-disjoint paths, and for every $2 \leq i \leq k$, the endpoints of $P_i$ are $v_i$ and $v_{2k-i+2}$.  
 
Finally, we establish (S3). Let $W = (\{v_1, \ldots, v_{2k}\} \cup V(P'_2) \cup \ldots \cup V(P'_k)) \setminus \{v'_1, v'_{k+1}\}$, and let $G_2 = G[V(G) \setminus W]$. Note that $N_G(\textrm{SMALL}) \setminus U \subseteq A \subseteq V(G) \setminus W$. It thus follows by Property (d) from Lemma~\ref{lem::OurCycleLemma} and by Theorem~\ref{th::MinDeg3} that $G_2$ satisfies the assertion of Lemma~\ref{lem::NewHamCon} and thus contains a Hamilton path whose endpoints are $v'_1$ and $v'_{k+1}$. Adding (if needed) the edges $v_1 v'_1$ and $v_{k+1} v'_{k+1}$ yields the required Hamilton path whose endpoints are $v_1$ and $v_{k+1}$.
\end{proof}

As noted at the beginning of this section, we end it with a short explanation of how one could adjust our proof of Theorem~\ref{th::CycleSpaceGnp} to handle larger values of $p$. First, the assertion of this theorem holds for $p \geq C \log n/n$, where $C$ is a sufficiently large constant, by a result from~\cite{CNP}. Hence, given an appropriate value of $C$, it suffices to consider the range $2 \log n/n \leq p \leq C \log n/n$. Standard calculations show that for such values of $p$ a.a.s. $\delta(G) = \Omega(\log n)$, $\Delta(G) = O(\log n)$ (the constant depends on $C$), and properties (P4)--(P6) remain valid; Property (P3) is not needed. Additionally, Lemma~\ref{lem::robustEdgeBetweenSets} still holds, since its proof only relies on Property (P6) whose proof does not require  an upper bound on $p$. Then, a variant of Lemma~\ref{lem::DiameterR} without any exceptional small vertices holds and provides short paths between any given pair of vertices even after a few vertices are deleted. Lemma~2.3 from~\cite{CNP} then yields an even cycle of length $O(\log n/\log\log n)$ containing an odd number of edges of $R$.
Partitioning the remaining vertices into two almost equal parts using 
Lemma~\ref{lem::LLLsplit}, we apply Theorem~\ref{th::PathsInExpanders} to construct the connecting paths in one part, with parameters $n_0 = n/\log n$ and $d = \sqrt{\log n}/2$. All vertices of the other part remain available and every vertex of the residual graph retains $\Omega(\log n)$ neighbours there. Lemma~\ref{lem::subgraphExpander} and Property (P6) imply that this residual graph is a sufficiently strong expander, and thus Theorem~\ref{th::HamConCexpander} yields the required Hamilton path.

\section*{Acknowledgements} We thank the anonymous referees for their helpful comments.



\begin{thebibliography}{99}

\bibitem{ALW}
B. Alspach, S. C. Locke, and D. Witte, The Hamilton spaces of Cayley graphs on abelian groups, \emph{Discrete Mathematics} 82 (2) (1990), 113--126.

\bibitem{BK}
J. D. Baron and J. Kahn, On the cycle space of a random graph, \emph{Random Structures and Algorithms} 54(1) (2019), 39--68.

\bibitem{BFHK}
S. Ben-Shimon, A. Ferber, D. Hefetz, and M. Krivelevich, Hitting time results for Maker-Breaker games,
\emph{Random Structures and Algorithms} 41 (2012), 23--46.

\bibitem{Bollobas}
B. Bollob\'as, The evolution of sparse graphs, in \emph{Graph Theory and Combinatorics}
(Cambridge, 1983), Academic Press, London, (1984), 35--57.

\bibitem{BL}
J. A. Bondy and L. Lov\'asz, Cycles through specified vertices of a graph, \emph{Combinatorica} 1 (1981), 117--140.

\bibitem{CNP}
M. Christoph, R. Nenadov, and K. Petrova, The Hamilton space of pseudorandom graphs, \emph{Journal of Combinatorial Theory, Series B} 176 (2026), 254--267.

\bibitem{DHJ}
B. DeMarco, A. Hamm, and J. Kahn, On the triangle space of a random graph, \emph{Journal of Combinatorics} 4(2) (2013), 229--249.

\bibitem{DGMS}
N. Dragani\'c, S. Glock, D. Munha Correia, and B. Sudakov, Optimal Hamilton covers and linear arboricity for random graphs, \emph{Proceedings of the American Mathematical Society} 153 (2025), 921--935.

\bibitem{DKN}
N. Dragani\'c, M. Krivelevich, and R. Nenadov, Rolling backwards can move you forward: on embedding problems in sparse expanders, \emph{Transactions of the American Mathematical Society} 375 (7) (2022), 5195--5216.

\bibitem{DMMPS}
N. Dragani\'c, R. Montgomery, D. Munha Correia, A. Pokrovskiy, and B. Sudakov, Hamiltonicity of expanders: optimal bounds and applications, arXiv preprint arXiv:2402.06603v2, 2024.

\bibitem{ER}
P. Erd\H{o}s and A. R\'enyi, On the evolution of random graphs, \emph{Bull. Inst. Statist. Tokyo} 38 (1961), 343--347.

\bibitem{FKL}
A. Ferber, G. Kronenberg, and E. Long, Packing, counting and covering Hamilton cycles
in random directed graphs, \emph{Israel Journal of Mathematics} 220 (2017), 57--87.

\bibitem{FK}
A. Frieze and M. Karo\'nski, \textbf{Introduction to Random Graphs}, Cambridge University Press, 2015.

\bibitem{Hartman}
I. B.-A. Hartman, Long cycles generate the cycle space of a graph, \emph{European Journal of Combinatorics} 4 (1983), 237--246.

\bibitem{HKT}
D. Hefetz, M. Krivelevich and T. Szab\'o, Sharp threshold for the appearance of certain spanning trees in random graphs,
\emph{Random Structures and Algorithms} 41 (2012), 391--412.

\bibitem{HKLO}
D. Hefetz, D. K\"uhn, J. Lapinskas, and D. Osthus, Optimal covers with Hamilton cycles in random graphs,
\emph{Combinatorica} 34 (2014), 573--596.

\bibitem{Heinig1}
P. Heinig, On prisms, M\"obius ladders and the cycle space of dense graphs, \emph{European Journal of
Combinatorics} 36 (2014), 503--530.

\bibitem{Heinig2}
P. Heinig, When Hamilton circuits generate the cycle space of a random graph, arXiv preprint arXiv:1303.0026, 2013.

\bibitem{JLR}
S. Janson, T. \L uczak, and A. Ruci\'nski, \textbf{Random graphs}, Wiley-Interscience Series in Discrete Mathematics and Optimization, Wiley-Interscience, New York, 2000.

\bibitem{KKO}
F. Knox, D. K\"uhn, and D. Osthus, Edge-disjoint Hamilton cycles in random graphs, \emph{Random Structures and Algorithms} 46 (2015), 397--445.

\bibitem{KSz}
J. Koml\'os and E. Szemer\'edi, Limit distributions for the existence of Hamilton circuits
in a random graph, \emph{Discrete Mathematics} 43 (1983), 55--63.

\bibitem{Korshunov}
A. D. Korshunov, Solution of a problem of Erd\H{o}s and R\'enyi on Hamilton cycles in non-oriented graphs, \emph{Soviet Math. Dokl.} 17 (1976), 760--764.

\bibitem{KSa}
M. Krivelevich and W. Samotij, Optimal packings of Hamilton cycles in sparse random graphs, \emph{SIAM Journal on Discrete Mathematics} 26 (2012), 964--982.

\bibitem{LS}
C. Lee and B. Sudakov, Dirac's theorem for random graphs, \emph{Random Structures and Algorithms} 41 (2012), 293--305.

\bibitem{Locke1}
S. C. Locke, A basis for the cycle space of a 2-connected graph, \emph{European Journal of Combinatorics} 6 (1985), 253--256.

\bibitem{Locke2}
S. C. Locke, A basis for the cycle space of a 3-connected graph, \emph{Annals of Discrete Mathematics} 27 (1985), 381--397.

\bibitem{Mont}
R. Montgomery, Hamiltonicity in random graphs is born resilient, \emph{Journal of Combinatorial Theory, Series B} 139 (2019), 316--341.

\bibitem{Posa}
L. P\'osa, Hamiltonian circuits in random graphs, \emph{Discrete Mathematics} 14 (1976), 359--364.

\end{thebibliography}
\end{document}